\documentclass[a4paper,12pt]{amsart}
\usepackage[latin1]{inputenc}
\usepackage{enumerate,amsfonts,amsmath,amsthm,amssymb,latexsym,verbatim,amscd,mathrsfs,hyperref}
\usepackage[usenames]{color}
\usepackage[all]{xy}

\usepackage[left=3cm,right=3cm,top=3cm,bottom=3cm]{geometry}

\numberwithin{equation}{section}
\theoremstyle{definition}
\newtheorem{thm}{Theorem}
\newtheorem{prop}[thm]{Proposition}
\newtheorem{cor}[thm]{Corollary}
\newtheorem{lem}[thm]{Lemma}
\newtheorem{defn}[thm]{Definition}
\newtheorem{rem}[thm]{Remark}


 
\title{Identities and isomorphisms of finite-dimensional graded simple algebras}
\author{Angelo Bianchi}
\email{acbianchi@unifesp.br}
\address{Instituto de Ci\^encia e Tecnologia, Universidade Federal de S\~ao Paulo, S\~ao Jos\'e dos Campos/SP, 12247-014, Brazil}
\thanks{A.B. is partially supported by FAPESP grant 2014/09310-5}

\author{Diogo Diniz}
\email{diogo@mat.ufcg.edu.br}
\address{Unidade Acad\^emica de Matem\'atica, Universidade Federal de Campina Grande, Campina Grande/PB, 58429-970, Brazil}
\thanks{D. D. is partially supported by  CNPq grants ~303822/2016-3 and ~406401/2016-0}

\begin{document}
\begin{abstract} 
Let $\mathbb F$ be an algebraically closed field, $G$ be an abelian group, and let $A$ and $B$ be arbitrary finite-dimensional $G$-graded simple  algebras over $\mathbb F$. We prove that $A$ and $B$ are isomorphic as graded algebras if, and only if, they satisfy the same graded polynomial identities.
\end{abstract}
\maketitle

\section*{Introduction}
The study of isomorphisms between algebras with same polynomial identities has intensely developed in the last five decades. Clearly, two non-isomorphic algebras $A$ and $B$ may have the same identities. But, on the other hand, even if both $A$ and $B$ are finite-dimensional simple algebras with same identities, it does not follow that $A$ and $B$ are isomorphic. Even more rigid structures such as Lie algebras can have the same identities not being isomorphic.

However, under the hypothesis that the base field is algebraically closed, the situation is quite different: as a consequence of the famous Amitsur-Levitsky, two finite dimensional simple associative algebras $A$ and $B$ have the same identities if, and only if, $A$ and $B$ are isomorphic. 

The purpose of this article is to prove that finite-dimensional (associative and non-associative) simple $G$-graded algebras over an algebraically closed field $\mathbb F$ of arbitrary characteristic and an abelian group $G$ are determined, up to $G$-graded isomorphism, by their $G$-graded identities. Two similar results about graded simple associative algebras are known:
\begin{itemize}
\item $\mathbb F$ algebraically closed in characteristic zero and arbitrary group $G$, by Aljadeff and Haile in \cite{AH};

\item $\mathbb F$ algebraically closed and $G$ abelian such that the orders of all finite subgroups of $G$ are invertible in $\mathbb F$, by Koshlukov and Zaicev in \cite{KZ},
\end{itemize}

Analogous results were obtained for Lie algebras by Kushkulei and Razmyslov in \cite{KR}, for Jordan algebras by Drensky and Racine in \cite{DR}, and for general non-associative (and associative) algebras by Shestakov and Zaicev in \cite{SZ}.

\

More specifically, the main result of this paper is:

\begin{thm} \label{main} Let $\mathbb F$ be an algebraically closed field and let $G$ be an abelian group. Then, two finite dimensional $G$-graded simple (associative or non-associative) $\mathbb F$-algebras have the same graded identities if, and only if, they are isomorphic.
\end{thm}

The proof of Theorem \ref{main} is performed for associative and non-associative algebras in a uniform manner, inspired by the approach of \cite{SZ}. This theorem  recovers \cite{SZ} and \cite{KZ} as particular cases. Further, the hypothesis on the group $G$ is more general than the hypothesis in \cite{KZ}. In fact, in \cite{KZ} the hypothesis that the orders of all finite subgroups of $G$ must be invertible in $\mathbb F$, comes from the description of all graded simple associative algebras, according to \cite[Theorems 2 and 3]{BZS}. We note that no classification is used here. The commutativity of the group $G$ is necessary to verify that the multiplication algebra $M(\mathcal{U})$ of a $G$-graded algebra $\mathcal U$ is a homogeneous subalgebra of $\mathrm{End}_{\mathbb{F}} \mathcal U^{gr}$. In this case $M(\mathcal{U})$ has a grading that makes $\mathcal U$ a graded $M(\mathcal{U})$-module. The hypothesis on the ground field will be crucial to show Lemma \ref{1} which is indispensable to get Theorem \ref{main}. 

In Section \ref{prereq} we establish basic definitions and in Section \ref{results} we discuss some intermediate results. We prove that a graded simple algebra $\mathcal U$ such that the multiplication algebra $M(\mathcal U)$ satisfies an ordinary polynomial identity is a finite dimensional graded vector space over the graded centroid $\Gamma(\mathcal U)$ (see Definition \ref{centr}). We prove, as a consequence of the Amitsur-Levitzki theorem that the dimension $\mathrm{dim}_{\Gamma(\mathcal U)}\ \mathcal U$ is determined by the graded identities of $\mathcal U$. Our main result, Theorem \ref{main}, is proved in Section \ref{proof}.

\

\noindent \textbf{Acknowledgment:} We wish to thank Dr. Ivan Shestakov, Dr. Thiago Castilho, and Dr. Claudemir Fidelis for helpful suggestions and discussions.

\section{Preliminaries} \label{prereq}
Let $G$ be a group, with unit element $e$, and $\mathbb{F}$ a field. An $\mathbb{F}$-algebra $\mathcal U$ is said to be $G$-{\bf graded} if there is a vector space decomposition $\mathcal U=\bigoplus_{g\in G} {\mathcal U}_g$ and ${\mathcal U}_g {\mathcal U}_h\subseteq {\mathcal U}_{gh}$ for all $g,h\in G$.  We denote by  ${\rm supp}(\mathcal U)=\{g\in G\mid \mathcal U_g\ne 0\}$ the \textbf{support of $G$}. A non-zero element $u\in {\mathcal U}$ is called {\bf homogeneous of degree} $g$ if $a\in A_g$ and, in this case, we write $\deg_G a=g$. A subspace $V\subseteq {\mathcal U}$ is called {\bf homogeneous} when $V=\bigoplus_{g\in \Gamma}(V\cap {\mathcal U}_g)$. We call ${\mathcal U}$ {\bf graded simple} if $\mathcal U^2\neq 0$  and $\mathcal U$ has no nontrivial homogeneous ideals. If an associative graded algebra $A$ has a unit and every non-zero homogeneous element is invertible we say that $A$ is a \textbf{graded division algebra}. Notice that any simple algebra with an arbitrary grading is graded simple. A $G$-graded simple algebra may not be simple as an ungraded algebra, for instance the group algebra $A=\mathbb{F}[G]$ of a finite group $G$ is graded simple with its canonical $G$-grading, however it is not simple as an ungraded algebra. We say that two $G$-graded algebras ${\mathcal U}=\oplus_{g\in G} {\mathcal U}_g$ and ${\mathcal V} =\oplus_{g\in G} {\mathcal V}_g$ are \textbf{isomorphic} as graded algebras if there exists an isomorphism of algebras $\psi: {\mathcal U} \to {\mathcal V}$ such that $\psi({\mathcal U}_g) = {\mathcal V}_g$ for all $g\in G$. Homomorphisms of graded algebras are defined analogously.

We recall the concept of graded and ordinary polynomial identities: consider an absolutely free algebra ${\mathbb F}\{X\}$ with a free generating set $X=\cup_{g\in G} X_g$, where $\{X_g\mid g\in G\}$ is a family of disjoint countable sets. One can define a $G$-grading on ${\mathbb F}\{X\}$ by setting $\deg_G x= g$ when $x\in X_g$ and extend this to ${\mathbb F}\{X\}$ in the natural way. A polynomial  $f(x_1,\dots,x_n)$ in variables $x_1\in X_{g_1},\dots, x_n\in X_{g_n}$ is called a {\bf graded polynomial identity} for ${\mathcal U}$ if $f(u_1,\dots, u_n) = 0$ for any $u_1 \in {\mathcal U}_{g_1},\dots, u_{n}\in {\mathcal U}_{g_n}$. We say that a polynomial $f(x_1,\dots,x_n)$ in the absolutely free algebra ${\mathbb F}\{X\}$, freely generated by the countable set $X$, is an {\bf ordinary polynomial identity} of ${\mathcal U}$ if $f(u_1,\dots,u_n)=0$ for any $u_1,\dots,u_n\in {\mathcal U}$. If $A$ is an associative algebra then we consider the ordinary polynomial identities for $A$ as polynomials in the free associative algebra $\mathbb{F}\langle X \rangle$.

Let $R$ be an associative $\mathbb{F}$-algebra with a $G$-grading. A $G$-grading on a left $R$-module $V$ is a vector space decomposition $V=\oplus_{g\in G} V_g$ such that $R_gV_h\subseteq V_{gh}$, for every $g,h \in G$. The module $V$ is {\bf graded simple} (or graded irreducible) if $RV\neq 0$ and the only homogeneous submodules of $V$ are $0$ and $V$. A {\bf graded vector space} $V$ is a graded module over a graded division algebra $R$.

If $\mathcal U$ is a $G$-graded $\mathbb{F}$-algebra and $\widetilde{\mathbb{F}}$ is an extension of $\mathbb{F}$ then in the algebra $\widetilde{\mathcal U}=\widetilde{\mathbb{F}}\otimes_{\mathbb{F}} \mathcal U$ we consider the $G$-grading  $\widetilde{\mathcal U}=\oplus_{g\in G}\widetilde{\mathcal U}_g$, where $\widetilde{\mathcal U}_g=\widetilde{\mathbb{F}}\otimes_{\mathbb{F}} \mathcal U_g$.

\section{Intermediate results} 

The goal of this section is to discuss some results about graded algebras, multiplication algebras of a given algebra, graded algebras with ordinary polynomial identities, and some considerations on the graded centroid of a graded algebra.

\subsection{On multiplication algebras}\label{results}

The multiplication algebra of an $\mathbb{F}$-algebra $\mathcal{U}$, denoted $M(\mathcal{U})$, is the subalgebra of $\mathrm{End}_{\mathbb{F}}\ \mathcal{U}$ generated by the linear transformations corresponding to left and right multiplications by elements of $\mathcal{U}$. Then $\mathcal{U}$ is an $M(\mathcal{U})$-module, the ring of endomorphisms of this module is the centroid of $\mathcal{U}$.

If $\mathcal{U}$ is graded by an abelian group $G$, then the left and right multiplications by homogeneous elements of $\mathcal U$ are homogeneous linear maps, therefore $M(\mathcal U)$ is a homogeneous subalgebra of $\left(\mathrm{End}_{\mathbb{F}} \mathcal U \right)^{gr}$. The algebra $\mathcal U$ becomes a graded $M(\mathcal U)$-module and the homogeneous submodules of this module are the homogeneous ideals of $\mathcal U$. 

\begin{defn}\label{centr}
Let $\mathcal U$ be an $\mathbb{F}$-algebra graded by an abelian group $G$ and let $M(\mathcal U)$ be the multiplication algebra of $\mathcal U$. The graded centroid of $\mathcal U$ is the $\mathbb{F}$-algebra $\Gamma(\mathcal U)=\left(\mathrm{End}_{M(\mathcal U)} \mathcal U \right)^{gr}$.
\end{defn}

\begin{lem}[Graded Schur Lemma] \cite[Lemma 2.4]{EK} \label{grschur}
Let $R$ be a $G$-graded associative algebra. Suppose $V$ is a graded simple left
$R$-module. Let $D=\left(\mathrm{End}_{R} V \right)^{gr}$. Then, $D$ is a graded division algebra. \hfill \qedsymbol
\end{lem}

Note that $\Gamma(\mathcal U)$ is a subalgebra $\mathrm{End}_{M(\mathcal U)} \mathcal U$, which is the centroid of $\mathcal U$ as an ungraded algebra. It follows from \cite[Proposition 5.4.1]{J2}  that if $\mathcal U^2 =\mathcal U$ then $\mathrm{End}_{M(\mathcal U)} \mathcal U$, and a fortiori  $\Gamma (\mathcal U)$, is commutative.  Note that if $\mathcal U$ is graded simple then $\mathcal U^2$ is a non-zero homogeneous ideal, hence $\mathcal U^2=\mathcal U$. Together with the graded version of Schur's Lemma this implies that  $\Gamma(\mathcal U)$ is a commutative graded division ring if $\mathcal U$ is a graded simple algebra.  

\begin{rem}\label{algclosed}
Let $\mathbb{F}$ be an algebraically closed field. If $\mathcal U$ is a graded simple finite dimensional $\mathbb{F}$-algebra, then $\Gamma(\mathcal U)_e$ is a field. Note that $\mathbb{F}\subseteq \Gamma(\mathcal U)_e$ is a finite extension, therefore $\Gamma(\mathcal U)_e=\mathbb{F}$. 
\end{rem}

\begin{thm}[Graded Jacobson density Theorem] \label{dens}\cite[Theorem 2.5]{EK}
Let $R$ be an associative $G$-graded algebra. Suppose $V$ is a graded simple left
$R$-module and let $D = End_
R^{gr}(V )$. If $v_1,\dots, v_n \in V$ are homogeneous elements that
are linearly independent over $D$, then for any $w_1,\dots,  w_n \in V$ there exists $r \in R$
such that $rv_i = w_i$, $i = 1,\dots, n$.
\end{thm}

\begin{cor}\label{endv}
Let $R$ be an associative $G$-graded algebra. Suppose $V$ is a graded simple faithful left
$R$-module and let $D = End_
R^{gr}(V )$. Let $\psi:R\rightarrow \mathrm{End}_D^{gr}V$ be the homomorphism of graded algebras such that $\psi(r)(v)=rv$ for every $r\in R$ and every $v\in V$. If $R$ satisfies an ordinary polynomial identity, then $V$ is a finite dimensional graded vector space over $D$.  Moreover in this case $\mathrm{End}_D^{gr}V=\mathrm{End}_DV$ and $\psi$ is an isomorphism of graded algebras.
\end{cor}
\proof
Let $\psi:R\rightarrow \mathrm{End}_D^{gr}V$ be the homomorphism of graded algebras such that $\psi(r)(v)=rv$ for every $r\in R$ and every $v\in V$. This is an injective homomorphism  because $V$ is a faithful module. It follows from \cite[Corollary I.2.11]{NV} that, if $V$ is finite-dimensional over $D$, then $\mathrm{End}_D^{gr}V=\mathrm{End}_DV$. Moreover, Theorem \ref{dens} implies that this homomorphism is also surjective and the result follows. We assume that there exists a sequence $v_1,v_2,\dots $ of homogeneous elements of $V$ linearly independent over $D$. Then Theorem \ref{dens} implies that for every $n$ there exists $r_{ij}$, $1\leq,i,j \leq n$ such that $r_{ij} v_k=\delta_{jk} v_i$. Let $S$ be the subalgebra of $R$ generated by the elements $r_{ij}$, $1\leq i,j \leq n$. The $\mathbb{F}$-subspace $W$ of $V$ generated by the elements $v_1,\dots, v_n$ is an $S$-module. Therefore we obtain a homomorphism  $\varphi:S\rightarrow \mathrm{End}_{\mathbb{F}}\ W$ of algebras such that $\varphi(s)w=sw$, for every $s\in S$ and every $w\in W$. This homomorphism maps $r_{ij}$ to the linear map $e_{ij}$ given by $e_{ij}v_k=\delta_{jk} v_i$, therefore it is a surjective homomorphism. Since $\mathrm{End}_{\mathbb{F}} W$ is isomorphic as an ungraded algebra to $M_n(\mathbb{F})$ this implies that an ordinary polynomial identity $f$ for $R$ is also a polynomial identity for $M_n(\mathbb{F})$. Since $n$ is arbitrary we may take $n>d/2$, where $d$ is the degree of $f$. This is a contradiction because $M_n(\mathbb{F})$ satisfies no identity of degree less than $2n$.
\endproof

The next proposition is a graded analogue of a result due to Jacobson, in accordance with \cite[ Corollary of Theorem 3]{J}.

\begin{prop}\label{grjac}
Let $\mathcal U$ be a graded algebra such that $M(\mathcal U)$ satisfies an ordinary polynomial identity. The algebra $\mathcal U$ is a graded simple algebra if, and only if, $M(\mathcal U)$ is a graded simple algebra. Moreover, $\mathcal{U}$ has finite dimension as a vector space over $\Gamma(\mathcal U)$ and  $M(\mathcal U)=\mathrm{End}_{\Gamma({\mathcal U})}\mathcal U$.
\end{prop}
\proof
Let $R=M(\mathcal U)$.   If $\mathcal U$ is a graded simple algebra then it is a faithful graded irreducible $R$-module. Therefore it follows from Corollary \ref{endv} that $\mathcal{U}$ has finite dimension as a vector space over $\Gamma(\mathcal U)$ and that $M(\mathcal U)=\mathrm{End}_{\Gamma{(\mathcal U})}\mathcal U$ is a graded simple algebra. To prove the converse note that any unital graded simple algebra admits a faithful graded irreducible module. Let $V$ be a graded irreducible faithful module for $R$ and let $D = End_
R^{gr}(V )$. Corollary \ref{endv} implies that $V$ is finite-dimensional over $D$ and that $R$ is isomorphic to $\mathrm{End}_DV$. Therefore the graded left module $_R R$ is completely reducible. This implies that any graded left $R$-module is completely reducible. In particular, $\mathcal{U}$ is a completely reducible graded module.  If $\mathcal U$ is not a graded simple algebra then  there exist non-zero graded submodules $\mathcal U_1$  and $\mathcal U_2$ such that $\mathcal U = \mathcal U_1 \oplus \mathcal U_2$. Denote by $M_i$ the subalgebra of $M(\mathcal U)$ generated by left and right multiplications by elements of $\mathcal U_i$, $i=1,2$. Then $M_1$ and $M_2$ are homogeneous ideals of $M(\mathcal U)$ such that $M(\mathcal U )=M_1\oplus M_2$, this is a contradiction. Hence, $\mathcal U$ is a graded simple algebra. 
\endproof

The previous result implies that a graded simple algebra $\mathcal U$ with a multiplication algebra $M(\mathcal U)$ that satisfies an ordinary polynomial identity is a finite dimensional vector space over the commutative graded division ring $\Gamma (\mathcal U)$. We now take a brief detour from our main goal to prove in Corollary \ref{eqdim}, as a consequence of the Amitsur-Levitzki theorem, that the dimension $\mathrm{dim}_{\Gamma(\mathcal U)} \mathcal U$ is determined by the ordinary identities of $\mathcal U$. We remark that it is well known that the graded identities of an algebra determine the ordinary identities. More precisely, if two graded algebras satisfy the same graded identities then they satisfy the same ordinary identities.

\begin{prop}\label{idm}
If the algebras $\mathcal{U}$ and $\mathcal V$ satisfy the same  polynomial identities then $M(\mathcal U)$ and $M(\mathcal V)$ satisfy the same polynomial identities
\end{prop}
\proof
Let $f(x_1,\dots, x_n)$ be a polynomial identity for $M(\mathcal U)$ and let $\overline{x_1},\dots, \overline{x_n}$ be elements in $M(\mathcal V)$. The algebra $M(\mathcal V)$ is generated by the left and right multiplications by elements of $\mathcal V$, therefore there exist polynomials $g_1,\dots, g_n$ in the free associative algebra $F\langle X \rangle$ and elements $v_1,\dots, v_r, w_1,\dots, w_s \in \mathcal V$ such that $\overline{x_i}=g_i((v_1)_L,\dots, (v_r)_L,(w_1)_R,\dots, (w_s)_R)$, $i=1,\dots, n$. We have $$f(\overline{x_1},\dots, \overline{x_n})=g((v_1)_L,\dots, (v_r)_L,(w_1)_R,\dots, (w_s)_R),$$
where $g(x_1,\dots, x_r, x_{r+1},\dots, x_{r+s})=f(g_1,\dots, g_n)$. Let $g_0$ denote the element $$g((x_1)_L,\dots, (x_r)_L,(x_{r+1})_R,\dots, (x_{r+s})_R)$$ in the multiplication algebra $M(F\left\{X\right\})$ of the free algebra $F\left\{X\right\}$. Then $g_0\cdot x_{r+s+1}$ is a polynomial in $F\left\{ X \right\}$. Since $f$ is an identity for $M(\mathcal U)$ we conclude that $g_0\cdot x_{r+s+1}$ is a polynomial identity for $\mathcal U$. Therefore $g_0\cdot x_{r+s+1}$ is a polynomial identity for $\mathcal V$. This implies that $$g((v_1)_L,\dots, (v_r)_L,(w_1)_R,\dots, (w_s)_R)\cdot v=0,$$ for every $v\in \mathcal V$. Hence  $g((v_1)_L,\dots, (v_r)_L,(w_1)_R,\dots, (w_s)_R)=0$. Therefore $f$ is a polynomial identity for $M(\mathcal V)$. Analogously, every polynomial identity for $M(\mathcal V)$ is also a polynomial identity for $M(\mathcal U)$. 
\endproof

\begin{cor}\label{eqdim}
Let $\mathcal U$ and $\mathcal V$ be graded simple algebras such that $M(\mathcal U)$ and $M(\mathcal V)$ are ordinary p.i. algebras. If $\mathcal U$ and $\mathcal V$ satisfy the same ordinary polynomial identities, then $\mathrm{dim}_{\Gamma(\mathcal U)} \mathcal U= \mathrm{dim}_{\Gamma(\mathcal V)} \mathcal V$.
\end{cor}
\proof
Proposition \ref{grjac} implies $\mathcal U$ and $\mathcal V$ are finite dimensional over $\Gamma(\mathcal U)$ and $\Gamma(\mathcal V)$, respectively, and that  $M(\mathcal U)=\mathrm{End}_{\Gamma{(\mathcal U})}\mathcal U$ and $M(\mathcal V)=\mathrm{End}_{\Gamma(\mathcal V)}\mathcal V$. Let $n_{\mathcal U}=\mathrm{dim}_{\Gamma(\mathcal U)}\mathcal U$ and $n_{\mathcal V}=\mathrm{dim}_{\Gamma(\mathcal V)}\mathcal V$ . Note that $M(\mathcal U)$ and $M(V)$ are isomorphic as an ordinary algebras to the algebras $M_{n_{\mathcal U}}(\Gamma (\mathcal U))$, of $n_{\mathcal U}\times n_{\mathcal U}$ matrices over $\Gamma(U)$, and $M_{n_{\mathcal V}}(\Gamma (\mathcal V))$, of $n_{\mathcal V}\times n_{\mathcal V}$ matrices over $\Gamma(V)$, respectively.  As a consequence of Proposition \ref{idm}, we conclude that these algebras satisfy the same ordinary polynomial identities. Recall that $\Gamma(U)$ and $\Gamma(V)$ are commutative, hence the Amitsur-Levitzki Theorem implies that $n_{\mathcal U}=n_{\mathcal V}$.
\endproof

For an algebra $\mathcal{A}$, we denote by $Z(\mathcal{A})$ the center of A. We denote by $Z(\mathcal{A})_e$ the set of all elements of degree $e$ inside $Z(\mathcal{A})$. In the next proposition we adapt the proof of L. H. Rowen in \cite{R} that any ideal in an associative semiprime p.i. ring intersects the center non-trivially. The next proposition is a weaker version of the graded analogue of this result (see \cite{K}) that will be sufficient for our intended application.

\begin{prop}\label{gposner}
Let $\mathcal A$ be an associative $\mathbb{F}$-algebra graded by a group $G$ and assume that there exists an extension $\widetilde{\mathbb{F}}$ of the field $\mathbb{F}$ and a graded simple $\widetilde{\mathbb{F}}$-algebra $\widetilde{\mathcal A}$ such that:
\begin{enumerate}
\item[(i)] $\widetilde{\mathcal A}$ is isomorphic as a graded $\widetilde{\mathbb{F}}$-algebra to the ring of endomorphisms $\mathrm{End}_{D}V$ of a finite dimensional graded vector space over a commutative graded division $\widetilde{\mathbb{F}}$-algebra $D$;
\item[(ii)] $A$ is a homogeneous $\mathbb{F}$-subalgebra of $\widetilde{\mathcal A}$;
\item[(iii)] $\widetilde{\mathbb{F}}\mathcal A=\widetilde{\mathcal A}$.
\end{enumerate}
 If $Z(\mathcal A)_e$ is a field, then $\mathcal A$ is a graded simple $\mathbb{F}$-algebra that satisfies an ordinary polynomial identity.
\end{prop}
\proof
We prove that any homogeneous ideal intersects $Z(\mathcal A)_e$ non-trivially. If $Z(\mathcal A)_e$ is a field this clearly implies that $\mathcal A$ is graded simple. Let $I$ be a non-zero homogeneous ideal of the $\mathbb{F}$-algebra $\mathcal A$, then $\widetilde{\mathbb{F}}I$ is a non-zero homogeneous ideal in $\widetilde{\mathcal A}$. Since $\widetilde{\mathcal A}$ is a graded simple algebra we conclude that $\widetilde{\mathbb{F}}I=\widetilde{\mathcal A}$. There exist natural numbers $q, n_1,\dots, n_q$ such that the neutral component $\widetilde{\mathcal A}_e$ is isomorphic to $M_{n_1}(\mathbb{L})\oplus\cdots \oplus M_{n_q}(\mathbb{L})$, where $\mathbb{L}$ denotes the field $D_e$  (see the remarks before Theorem 1 in \cite{KZ}). Let $m=\mathrm{max}\{n_1,\dots, n_q\}$ and let $f(x_1,\dots, x_n)$ be a multihomogeneous central polynomial for $M_m(\mathbb{L})$ with coefficients in the prime field of $\mathbb{L}$ (for example the Formanek polynomial in \cite[Theorem 3.4]{DF}).  Then $f$ is a central polynomial for $\widetilde{\mathcal A}_e$. Denote $f(I_e)$ the $\mathbb{F}$-subspace of $\mathcal{A}$ generated by the elements $f(a_1,\dots, a_n)$ with $a_1,\dots, a_n \in I_e$. Note that $\widetilde{\mathcal A}_e=\widetilde{\mathbb{F}}I_e$, since $f$ is not an identity for $\widetilde{\mathcal A}_e$ we conclude that $f(I_e)\neq 0$. Clearly  $f(I_e)\subseteq Z(\mathcal A)_e\cap I$.  Since $\widetilde{\mathcal A}$ is isomorphic to $\mathrm{End}_D\ V$ and $D$ is commutative $\widetilde{\mathcal A}$ satisfies as an $\mathbb{F}$-algebra the same ordinary identities as $M_n(\mathbb{F})$, where $n=\mathrm{dim}_D\ V$.  Any polynomial identity for $\widetilde{\mathcal A}$ is also an  identity for $\mathcal A$.
\endproof

\begin{rem}\label{u}
Let $\mathbb{F}\subseteq \widetilde{\mathbb{F}}$ be an extension of fields. Let $\widetilde{\mathcal U}$ be an $\widetilde{\mathbb{F}}$-algebra and $\mathcal{U}$ a homogeneous $\mathbb{F}$-subalgebra of $\widetilde{\mathcal U}$  such that $\widetilde{\mathbb{F}}\mathcal U=\widetilde{\mathcal U}$. Denote $M$ the $\mathbb{F}$-subalgebra of $M(\widetilde{\mathcal U})$ generated by the left and right multiplications by elements of $\mathcal U$. Then $M$ is a homogeneous subalgebra of $M(\widetilde{\mathcal U})$ such that $\widetilde{\mathbb{F}}M=M(\widetilde{\mathcal U})$. Moreover the map $\psi:M\rightarrow M(\mathcal U)$ such that $\psi(m)(x)=m(x)$ for every $m\in M$ and every $x\in \mathcal U$ is an isomorphism of graded $\mathbb{F}$-algebras.
\end{rem}

\begin{lem}\label{simple}
Let $\mathcal U$ be a graded simple algebra such that $M(\mathcal U)$ satisfies an ordinary polynomial identity and let $\widetilde{\mathbb{F}}$ be an extension of the field $\Gamma(\mathcal U)_e$. The algebra $\widetilde{\mathcal U}=\widetilde{\mathbb{F}}\otimes_{\Gamma(\mathcal U)_e}\mathcal U$  is a graded simple algebra and $M(\widetilde{\mathcal U})$ satisfies an ordinary polynomial identity. Moreover the center of $M(\widetilde{\mathcal U})$ is isomorphic as a graded algebra to $\widetilde{\mathbb{F}}\otimes_{\Gamma(\mathcal U)_e}\Gamma(\mathcal U)$.
\end{lem}
\proof
Let $\Gamma(\mathcal U)_e=\mathbb{F}$ and denote $M$ the $\mathbb{F}$-subagebra of $M(\widetilde{\mathcal U})$ generated by left and right multiplication by elements of $1\otimes_{\mathbb{F}} \mathcal U$. It is clear that the map $\psi:\widetilde{\mathbb{F}}\otimes_{\mathbb{F}}M\rightarrow M(\widetilde{\mathcal U})$ such that $\psi(\lambda\otimes m)=\lambda m$ for every $\lambda \in \widetilde{\mathbb{F}}$ and every $m\in M$ is a homomorphism of graded algebras. This homomorphism is surjective because $\widetilde{\mathbb{F}}M=M(\widetilde{\mathcal U})$. An element in $\mathrm{ker}\ \psi$ may be written as $\lambda_1\otimes m_1+\cdots+ \lambda_k\otimes m_k$, where $\lambda_1,\dots, \lambda_k$ are elements in $\widetilde{\mathbb{F}}$ linearly independent over $\mathbb{F}$ and $m_1,\dots, m_k$ are elements of $M$. We have
\begin{equation}\label{t}
\lambda_1m_1+\cdots + \lambda_k m_k=0.
\end{equation}
Given $u\in \mathcal U$ there exist $u_1,\dots, u_k \in \mathcal U$ such that $m_i(1\otimes u)=1\otimes u_i$.  Then, (\ref{t}) implies that $$\lambda_1\otimes u_1+\cdots + \lambda_k \otimes u_k=0.$$ Therefore we conclude that $u_1,\dots, u_k=0$.  Since the element $u$ is arbitrary we conclude that $m_1=\dots= m_k=0$. Hence $\mathrm{ker}\ \psi=0$. Note that $M$ is isomorphic as a graded algebra to $M(\mathcal U)$, therefore we conclude that $\widetilde{\mathbb{F}}\otimes_{\mathbb{F}}M(\mathcal U)$ is isomorphic as a graded algebra to $M(\widetilde{\mathcal U})$. 
Proposition \ref{grjac} implies that $M(\mathcal U)=\mathrm{End}_{\Gamma(\mathcal U)} \mathcal U$ and that $\mathcal U$ is finite dimensional over $\Gamma (\mathcal U)$. Let $\{u_1,\dots, u_n\}$ be a basis for $\mathcal U$ over $\Gamma  (\mathcal U)$ and let ${\bf g}=(g_1,\dots, g_n)$, where $g_i$ is the degree of $u_i$, $i=1,\dots, n$. The algebra $\mathrm{End}_{\Gamma(\mathcal U)} \mathcal U$ is isomorphic as a graded algebra to $M_n(\mathbb{F})\otimes_{\mathbb{F}} \Gamma(\mathcal U)$, where $M_n(\mathbb{F})$  has the elementary grading induced by ${\bf g}$ (see \cite{KZ}). Therefore, $\widetilde{\mathbb{F}}\otimes_{\mathbb{F}}\mathrm{End}_{\Gamma(\mathcal U)} \mathcal U$ is isomorphic as a graded algebra to $M_n(\mathbb{F})\otimes_{\mathbb{F}}\left(\widetilde{\mathbb{F}}\otimes_{\mathbb{F}}\Gamma(\mathcal U)\right)$. Note that $\widetilde{\mathbb{F}}\otimes_{\mathbb{F}}\Gamma(\mathcal U)$ is a graded division algebra, therefore $M_n(\mathbb{F})\otimes_{\mathbb{F}}\left(\widetilde{\mathbb{F}}\otimes_{\mathbb{F}}\Gamma(\mathcal U)\right)$ is a graded simple algebra. Hence we conclude that $M(\widetilde{\mathcal U})$ is isomorphic as a graded algebra to $M_n(\mathbb{F})\otimes_{\mathbb{F}}\left(\widetilde{\mathbb{F}}\otimes_{\mathbb{F}}\Gamma(\mathcal U)\right)$. This implies that the center of $M(\widetilde{\mathcal U})$ is isomorphic as a graded algebra to $\left(\widetilde{\mathbb{F}}\otimes_{\mathbb{F}}\Gamma(\mathcal U)\right)$. Moreover $M(\widetilde{\mathcal U})$ is a graded simple algebra that satisfies an ordinary polynomial identity, therefore Proposition \ref{grjac} implies that $\widetilde{\mathcal U}$ is graded simple. 
\endproof

In the next proposition we establish the graded analogue of the multiplication algebra properties that were obtained for ordinary algebras in \cite{PS}.

\begin{prop}\label{simpleloc}
Let $\mathcal U$ be a graded simple algebra of dimension $n$ over $\mathbb{F}=\Gamma(\mathcal U)_e$ and for $g\in G$ let $\{u_1^{g},\dots, u_{n_g}^{g}\}$ be a basis for $\mathcal{U}_g$ as a vector space over $\mathbb{F}$. Moreover let $\widetilde{\mathbb{F}}$ be an extension of the field of fractions of the polynomial algebra $\mathbb{F}[T]$, where $T=\{t_{i,j}^g\mid 1\leq i,j \leq \mathrm{dim}\ \mathcal{U}_g, g\in \mathrm{supp}\ \mathcal {U}\}$. Let $\widetilde{\mathcal U}=\widetilde{\mathbb F}\otimes_{\mathbb F}\mathcal U$ and denote $\mathcal{F}_{\mathcal  U}$ the $\mathbb{F}$-subalgebra of $\widetilde{\mathcal U}$ generated by the generic elements $y_i^{g}=t_{i1}^{g}\otimes u_1^{g}+\cdots +t_{in_g}^{g}\otimes u_{n_g}^{g}$, $i=1,\dots,\mathrm{dim}\ \mathcal{U}_g$, $g\in \mathrm{supp}\ \mathcal U$. 
\begin{enumerate}
\item[(i)] The neutral component $Z_{\mathcal  U}=Z\left(M(\mathcal{F}_{\mathcal U})\right)_e$ of the center of $M(\mathcal{F}_{\mathcal U})$  is isomorphic to an $\mathbb{F}$-subalgebra of $\widetilde{\mathbb{F}}$. 
\item[(ii)] The central localization $S_{\mathcal U}^{-1}\mathcal{F}_{\mathcal U}$,  where $S_{\mathcal U}=Z(M(\mathcal{F}_{\mathcal U}))_e\setminus \{0\}$,  is a graded simple algebra such that $M(S_{\mathcal U}^{-1}\mathcal{F}_{\mathcal U}$) satisfies an ordinary polynomial identity. 
\item[(iii)] $\Gamma(S_{\mathcal U}^{-1}\mathcal{F}_{\mathcal U})_e$ is isomorphic to the field of fractions of $Z_{\mathcal U}$.	
\item[(iv)]  The map $\widetilde{\mathbb F}\otimes_{K_{\mathcal U}}S_{\mathcal U}^{-1}\mathcal{F}_{\mathcal U}\rightarrow \widetilde{\mathcal U}$, where $K_{\mathcal U}=\Gamma(S_{\mathcal U}^{-1}\mathcal{F}_{\mathcal U})_e$,  given by $\widetilde{\lambda}\otimes x \mapsto \widetilde{\lambda}x$, for every $\widetilde{\lambda}\in \widetilde{\mathbb{F}}$ and every $x\in S_{\mathcal U}^{-1}\mathcal{F}_{\mathcal U}$, is an isomorphism of graded algebras.
\end{enumerate}
\end{prop}
\proof
Note that $\{1\otimes u_1^{g},\cdots, 1\otimes u_{n_g}^{g}\}$ is a basis for $\widetilde{\mathcal U}_g=1\otimes \mathcal U_g$ over $\widetilde{\mathbb F}$. The coefficients of the elements $y_1^{g},\dots, y_{n_g}^{g}$ relative to this basis form the matrix $(t_{ij}^g)$ with non-zero determinant. Therefore these elements are linearly independent over $\widetilde{\mathbb F}$. This implies that $\widetilde{\mathbb{F}} \mathcal{F}_{\mathcal U}=\widetilde{\mathcal U}$. Let $M$ be the $\mathbb{F}$-subalgebra of $M(\widetilde{\mathcal U})$ generated by left and right multiplication by elements of $\mathcal{F}_{\mathcal U}$. Remark \ref{u} implies that there exists an isomorphism $\psi:M(\mathcal{F}_{\mathcal U})\rightarrow M$ of graded $\mathbb{F}$-algebras, this isomorphism maps $Z_{\mathcal  U}$ onto $Z(M)_e$. The equality $\widetilde{\mathbb F}M=M(\widetilde{\mathcal U})$ implies that $Z(M)_e$ is an $\mathbb F$-subalgebra of $Z(M(\widetilde{\mathcal U}))_e$. Lemma \ref{simple} implies that $Z(M(\widetilde{\mathcal U}))_e=\widetilde{\mathbb{F}}$. Therefore [(i)] holds.

We identify $Z_{\mathcal  U}$ with $\psi(Z_{\mathcal  U})$, therefore $S_{\mathcal U}^{-1}\subseteq \widetilde{\mathbb{F}}^{\times}$, and $S_{\mathcal U}^{-1}M(\mathcal{F}_{\mathcal U})$ with the $\mathbb{F}$-subalgebra $\{s^{-1}m\mid s\in S_{\mathcal U}, m\in M\}$ of $M(\widetilde{\mathcal U})$. Since $\widetilde{\mathbb{F}}M=M(\widetilde{\mathcal U})$ we conclude that $\widetilde{\mathbb{F}}\left(S_{\mathcal U}^{-1}M(\mathcal{F}_{\mathcal U})\right)=M(\widetilde{\mathcal U})$. Lemma \ref{simple} implies that $M(\widetilde{\mathcal U})$ satisfies an ordinary polynomial identity. The neutral component of the center of $S_{\mathcal U}^{-1}M(\mathcal{F}_{\mathcal U})$ is the field of fractions of $Z_{\mathcal U}$, therefore Proposition \ref{gposner} implies that $S_{\mathcal U}^{-1}M(\mathcal{F}_{\mathcal U})$ is a graded simple algebra that satisfies an ordinary polynomial identity.
The map $\varphi:S_{\mathcal U}^{-1}M(\mathcal{F}_{\mathcal U})\rightarrow M(S_{\mathcal U}^{-1}\mathcal{F}_{\mathcal U})$ such that $\varphi(s^{-1}m)(t^{-1}x)=(st)^{-1}m(x)$, for every $s,t \in S_{\mathcal U}$, $m\in M(\mathcal{F}_{\mathcal U})$, $x\in \mathcal U$ is an isomorphism of graded algebras. Hence we conclude that $M(S_{\mathcal U}^{-1}\mathcal{F}_{\mathcal U})$ is a graded simple algebra that satisfies an ordinary polynomial identity.  Proposition \ref{grjac} implies that $S_{\mathcal U}^{-1}\mathcal{F}_{\mathcal U}$ is a graded simple algebra. This proves [(ii)].

It follows from Proposition \ref{grjac} that for any graded simple algebra, such that the multiplication algebra is ordinary p.i., the centroid is isomorphic as a graded algebra to the center of the multiplication algebra. Therefore $\Gamma(S_{\mathcal U}^{-1}\mathcal{F}_{\mathcal U})_e$ is isomorphic to the neutral component of the center of $S_{\mathcal U}^{-1}M(\mathcal{F}_{\mathcal U})$, which is the field of fractions of $Z_{\mathcal U}$. Hence we conclude [(iii)].

It follows from [(ii)] that $S_{\mathcal U}^{-1}\mathcal{F}_{\mathcal U}$ is a graded simple algebra such that $M(S_{\mathcal U}^{-1}\mathcal{F}_{\mathcal U})$ satisfies an ordinary polynomial identity. Lemma \ref{simple} implies that $\widetilde{\mathbb F}\otimes_{K_{\mathcal U}}S_{\mathcal U}^{-1}\mathcal{F}_{\mathcal U}$ is a graded simple algebra.  The map given by $\widetilde{\lambda}\otimes x \mapsto \widetilde{\lambda}x$, for $\widetilde{\lambda}\in \widetilde{\mathbb{F}}$ and $x\in S_{\mathcal U}^{-1}\mathcal{F}_{\mathcal U}$, is a surjective homomorphism of graded algebras from $\widetilde{\mathbb F}\otimes_{K_{\mathcal U}}S_{\mathcal U}^{-1}\mathcal{F}_{\mathcal U}$ onto $\widetilde{\mathcal U}$. Since the domain is graded simple this homomorphism is also injective. Therefore [(iv)] holds.
\endproof

\begin{rem}\label{galg}
Let $\mathcal U$ be a finite dimensional graded algebra. The algebra $\mathcal{F}_{\mathcal  U}$ in the previous proposition is called the algebra of generic elements of $\mathcal U$. If $\mathcal V$ is a finite dimensional graded algebra that satisfies the same graded identities as $\mathcal U$ then the algebra $\mathcal{F}_{\mathcal  V}$ of generic elements of $\mathcal V$ is isomorphic as a graded $\mathbb{F}$-algebra to $\mathcal{F}_{\mathcal  U}$.
\end{rem}

\subsection{A technical result on isomorphisms of certain scalar extensions}

Next result is crucial to prove Theorem \ref{main} and requires the algebraically closedness of the base field. It is based on Shestakov and Zaicev's ideas (for non-graded algebras) from a private communication.

\begin{lem} \label{1}
Let $\mathcal U$ and $\mathcal V$ be finite dimensional $G$-graded algebras over an algebraically closed field $\mathbb F$. Suppose that for a certain extension  $\mathbb F\subseteq \widetilde{\mathbb F}$ and for a certain $\mathbb F$-automorphism $\sigma$ of the field $\widetilde{\mathbb F}$ there exists a $G$-graded isomorphism of $\mathbb F$-algebras $\varphi: \widetilde{\mathbb F}\otimes_\mathbb F \mathcal U \to \widetilde{\mathbb F}\otimes_\mathbb F \mathcal V$ such that $\varphi(\alpha \otimes u)=\sigma(\alpha)\varphi(1\otimes u)$ for any $\alpha \in \widetilde{\mathbb F}$ and $u\in\mathcal U$. Then, $\mathcal U$ and $\mathcal V$ are isomorphic as graded algebras. 
\end{lem}

\proof If $\{u_1,\dots, u_n\}$ is an $\mathbb{F}$-basis for $\mathcal U$ then $\{\varphi(1\otimes u_1),\dots, \varphi(1\otimes u_n)\}$ is an $\widetilde{\mathbb F}$-basis for $\widetilde{\mathbb F}\otimes_\mathbb F \mathcal V$. This implies that $\mathrm{dim}_{\mathbb F}\mathcal U=\mathrm{dim}_{\mathbb F}\mathcal V$. Let $\{v_1,\dots, v_n\}$ be an $\mathbb{F}$-basis for $\mathcal V$. Then there exists $\alpha_{ij}\in \widetilde{\mathbb{F}}$, $1\leq i,j \leq n$, such that $\varphi(1\otimes u_i)=\sum_j\alpha_{ij}\otimes v_j$. Note that $\alpha=\mathrm{det}\ (\alpha_{ij})$ is a non-zero element of $\widetilde{\mathbb{F}}$.
Let $\mathbb E$ be the $\mathbb F$-subalgebra of $\widetilde{\mathbb{F}}$ generated by $\mathcal S=\{\alpha_{ij}, \alpha^{-1}\}$. Then $\mathbb E$ is a commutative $\mathbb F$-algebra generated by $n^2+1$ elements, therefore it is isomorphic to $\frac{\mathbb{F}[x_1,\dots, x_{n^2+1}]}{I}$, for an adequate prime ideal $I$ of $\mathbb{F}[x_1,\dots, x_{n^2+1}]$. Now, since $\mathbb F$ is algebraically closed field, the Hilbert--Nullstellensatz theorem implies that there exists $(\lambda_1,\dots,\lambda_{n^2+1} )\in \mathbb{F}^{n^2+1}$ such that $f(\lambda_1,\dots,\lambda_{n^2+1} )=0$ for every $f\in I$. Now let $\phi:\mathbb{F}[x_1,\dots, x_{n^2+1}]\rightarrow \mathbb{F}$ be the homomorphism such that $\phi(f)=f(\lambda_1,\dots,\lambda_{n^2+1} )$. Since $I\subseteq \mathrm{ker}\ \phi$ and $\mathbb{E}\cong \frac{\mathbb{F}[x_1,\dots, x_{n^2+1}]}{I}$ we obtain an homomorphism $\psi:\mathbb E\rightarrow \mathbb{F}$. This homomorphism can be uniquely extended to an $\mathbb F$-homomorphism $\widetilde\psi: \mathbb E\otimes_\mathbb F \mathcal V\to \mathcal V$ such that $\widetilde\psi(e\otimes v)=\psi(e)v$ for all $e\in\mathbb E$, $v\in \mathcal V$. Note that $\widetilde{\psi}$ is a homomorphism of graded $\mathbb{F}$-algebras. We identify $\mathcal U$ with the $\mathbb{F}$-subalgebra $1\otimes_\mathbb{F} \mathcal U$ of $\widetilde{\mathbb{F}}\otimes_\mathbb{F} \mathcal U$. Then $\varphi$ maps $\mathcal U$ onto $\mathbb{E}\otimes_{\mathbb{F}}\mathcal V$, moreover $u\mapsto \varphi(u)$ is a homomorphism of graded algebras. We compose this with $\widetilde\psi$ to obtain a homomorphism $\Theta:\mathcal U\rightarrow \mathcal V$ of graded $\mathbb{F}$-algebras. Note that $\Theta(u_i)=\sum_j\psi(\alpha_{ij})v_j$. Since $ \psi(\mathrm{det}\ (\alpha_{ij}))\psi(\alpha^{-1})=\mathrm{det}\ (\psi(\alpha_{ij}))\psi(\alpha^{-1})=1$ we conclude that $\mathrm{det}\ (\psi(\alpha_{ij}))\ne 0$ and, hence, $\Theta$ is an isomorphism.
\endproof

\section{Proof of Theorem \ref{main} } \label{proof}
Let $G$ be an abelian group and let $\mathbb{F}$ be an algebraically closed field. Let $\mathcal{U}$ and $\mathcal{V}$ be two finite dimensional algebras with a grading by the group $G$ such that $\mathcal{U}$ and $\mathcal{V}$ are graded simple algebras that satisfy the same graded polynomial identities. Remark \ref{algclosed} implies that $\Gamma(\mathcal U)_e=\mathbb{F}=\Gamma(\mathcal V)_e$. It is clear that $\mathrm{supp}\ \mathcal U=\mathrm{supp}\ \mathcal V$. We denote $\widetilde{\mathcal U}=\widetilde{\mathbb F}\otimes_{\mathbb F}\mathcal U$ and $\widetilde{\mathcal V}=\widetilde{\mathbb F}\otimes_{\mathbb F}\mathcal V$, where $\widetilde{\mathbb{F}}$ 
is the field of fractions of the polynomial algebra $\mathbb{F}[T]$ where $T=\{t_{i,j}^g\mid  1\leq i, j \leq m_g, g\in \mathrm{supp}\ \mathcal U\}$, here $m_g=\mathrm{max}\{\mathrm{dim}\ \mathcal{U}_g, \mathrm{dim}\mathcal{V}_g\}$.  Given $g\in \mathrm{supp}\ \mathcal {U}$ let $\{u_1^{g},\dots, u_{n_g}^{g}\}$ and $\{v_1^{g},\dots, v_{n^{\prime}_g}^{g}\}$ be basis for $\mathcal{U}_g$ and $\mathcal{V}_g$, respectively, as vector spaces over $\mathbb{F}$. Denote $\mathcal{F}_{\mathcal  U}$ the $\mathbb{F}$-subalgebra of $\widetilde{\mathcal U}$ generated by the generic elements $y_i^{g}=t_{i1}^{g}\otimes u_1^{g}+\cdots +t_{in_g}^{g}\otimes u_{n_g}^{g}$, $g\in \mathrm{supp}\ \mathcal U$, $i=1,\dots,n_g$. Analogously, let $\mathcal{F}_{\mathcal  V}$ be the $\mathbb{F}$-subalgebra of $\widetilde{\mathcal V}$ generated by the generic elements $z_i^{g}=t_{i1}^{g}\otimes v_1^{g}+\cdots +t_{in^{\prime}_g}^{g}\otimes v_{n^{\prime}_g}^{g}$, $g\in \mathrm{supp}\ \mathcal V$, $i=1,\dots,n^{\prime}_g$. Moreover let $Z_{\mathcal U}=Z(M(\mathcal{F}_{\mathcal U}))_e$, $S_{\mathcal U}=Z_{\mathcal U}\setminus \{0\}$ and  $Z_{\mathcal V}=Z(M(\mathcal{F}_{\mathcal V}))_e$, $S_{\mathcal V}=Z_{\mathcal V}\setminus \{0\}$. As noted in Remark \ref{galg} $\mathcal{F}_{\mathcal U}$ is isomorphic to $\mathcal{F}_{\mathcal V}$  as a graded $\mathbb{F}$-algebra. Let $\psi:\mathcal{F}_{\mathcal U}\rightarrow \mathcal{F}_{\mathcal V}$ be an isomorphism of graded algebras, it induces an isomorphism $\psi_M:M(\mathcal{F}_{\mathcal U})\rightarrow M(\mathcal{F}_{\mathcal V})$ of $\mathbb{F}$-graded algebras such that $\psi_M(m)(\psi(x))=\psi(mx)$, for every $m\in M(\mathcal{F}_{\mathcal U})$ and every $x\in \mathcal U$. This isomorphism maps $Z_{\mathcal U}$ onto $Z_{\mathcal {V}}$. Therefore we obtain an isomorphism of graded $\mathbb{F}$-algebras $\Psi:S_{\mathcal U}^{-1}\mathcal{F}_{\mathcal U}\rightarrow S_{\mathcal V}^{-1}\mathcal{F}_{\mathcal V}$ given by $\Psi(s^{-1}f)=\psi_M(s)^{-1}\psi(f)$. This isomorphism induces an isomorphism $K_{\mathcal U}\rightarrow K_{\mathcal V}$ of $\mathbb{F}$-algebras, where $K_{\mathcal U}=\Gamma(S_{\mathcal U}^{-1}\mathcal{F}_{\mathcal U})_e$ and $K_{\mathcal V}=\Gamma(S_{\mathcal U}^{-1}\mathcal{F}_{\mathcal U})_e$. It follows from [(i)] and [(iii)] in Proposition \ref{simpleloc} that $K_{\mathcal U}$ and $K_{\mathcal V}$ are isomorphic as $\mathbb{F}$-algebras to subfields $\widetilde{\mathbb{F}}_{\mathcal U}$ and $\widetilde{\mathbb{F}}_{\mathcal V}$ of $\widetilde{\mathbb{F}}$, respectively.
Then we obtain an isomorphism $\varsigma:\widetilde{\mathbb{F}}_{\mathcal U}\rightarrow \widetilde{\mathbb{F}}_{\mathcal V}$  of $\mathbb F$-algebras. Note that this implies the equality ${\rm trdeg}_\mathbb F \widetilde{\mathbb{F}}_{\mathcal U}={\rm trdeg}_\mathbb F \widetilde{\mathbb{F}}_{\mathcal V}$.  Let $\overline{\mathbb F}$ be the algebraic closure of $\widetilde{\mathbb F}$. Then, ${\rm trdeg}_\mathbb F \widetilde{\mathbb F} = {\rm trdeg}_\mathbb F \overline{\mathbb F}={\rm trdeg}_\mathbb F \widetilde{\mathbb{F}}_{\mathcal U} + {\rm trdeg}_{\widetilde{\mathbb{F}}_{\mathcal U}}\overline{\mathbb F}$ and also ${\rm trdeg}_\mathbb F \widetilde{\mathbb F} = {\rm trdeg}_\mathbb F \overline{\mathbb F}={\rm trdeg}_\mathbb F \widetilde{\mathbb{F}}_{\mathcal V} + {\rm trdeg}_{\widetilde{\mathbb{F}}_{\mathcal V}}\overline{\mathbb F}.$ From this we conclude that ${\rm trdeg}_{\widetilde{\mathbb{F}}_{\mathcal U}}\overline{\mathbb F}={\rm trdeg}_{\widetilde{\mathbb{F}}_{\mathcal V}}\overline{\mathbb F}.$
Now, since $\overline{\mathbb F}$ is algebraically closed, the previous equality  together with \cite[Chapter VII, Proposition 3.1]{L} imply that there exists automorphism $\sigma$ of $\overline{\mathbb F}$ that extends $\varsigma$. We then obtain an isomorphism of graded algebras $\overline{\mathbb F}\otimes_{K_{\mathcal U}}S^{-1}\mathcal{F}_{\mathcal U}\rightarrow \overline{\mathbb F}\otimes_{K_{\mathcal V}}S^{-1}\mathcal{F}_{\mathcal V}$  such that  $\widetilde{\lambda}\otimes x\mapsto \sigma(\widetilde{\lambda})\otimes \Psi(x)$. We compose with the isomorphisms $\overline{\mathbb F}\otimes_{K_{\mathcal U}}S^{-1}\mathcal{F}_{\mathcal U}\rightarrow \tilde{\mathcal U}$ and $\overline{\mathbb F}\otimes_{K_{\mathcal V}}S^{-1}\mathcal{F}_{\mathcal V}\rightarrow \tilde{\mathcal V}$ in [(iv)] of Proposition \ref{simpleloc} to obtain an isomorphism $\widetilde{\mathcal U}\rightarrow \widetilde{\mathcal V}$ of $\mathbb{F}$-algebras. This isomorphism satisfies the hypothesis of Lemma \ref{1} and the result follows. \hfill\qedsymbol

\vfill

\end{document}